\documentclass[reqno,centertags]{amsart}
\usepackage{amsmath,amsthm,amscd,amssymb}
\usepackage{latexsym}
\newcommand{\bbD}{{\mathbb{D}}}

\newcommand{\bbG}{{\mathbb{G}}}

\newcommand{\calQ}{{\mathcal Q}}
\newcommand{\calR}{{\mathcal R}}

\newcommand{\lb}{\label}
\newcommand{\f}{\frac}

\newcommand{\ol}{\overline}

\newcommand{\wti}{\widetilde  }

\newcommand{\book}{\text{\rm{book}}}

\newcommand{\s}{\text{\rm{s}}}
\newcommand{\sing}{\text{\rm{sing}}}

\newcommand{\bi}{\bibitem}

\newcommand{\beq}{\begin{equation}}
\newcommand{\eeq}{\end{equation}}
\newcommand{\ba}{\begin{align}}
\newcommand{\ea}{\end{align}}
\newcommand{\veps}{\varepsilon}

\newcounter{smalllist}
\newenvironment{SL}{\begin{list}{{\rm\roman{smalllist})}}{%
\setlength{\topsep}{0mm}\setlength{\parsep}{0mm}\setlength{\itemsep}{0mm}%
\setlength{\labelwidth}{2em}\setlength{\leftmargin}{2em}\usecounter{smalllist}%
}}{\end{list}}

\allowdisplaybreaks
\numberwithin{equation}{section}

\newtheorem{theorem}{Theorem}[section]

\newtheorem*{p2.1}{Proposition 2.1}
\newtheorem{proposition}[theorem]{Proposition}
\newtheorem{lemma}[theorem]{Lemma}

\theoremstyle{definition}
\newtheorem{example}[theorem]{Example}

\theoremstyle{remark}
\newtheorem*{remark}{Remark} 
\newtheorem*{remarks}{Remarks}

\theoremstyle{definition}
\newtheorem*{definition}{Definition}

\newcommand{\abs}[1]{\lvert#1\rvert}

\begin{document}
\title[Meromorphic Szeg\H{o} Functions]{Meromorphic Szeg\H{o} Functions 
and Asymptotic Series for Verblunsky Coefficients}
\author[B. Simon]{Barry Simon*}

\thanks{$^*$ Mathematics 253-37, California Institute of Technology, Pasadena, CA 91125.
E-mail: bsimon@caltech.edu. Supported in part by NSF grant DMS-0140592 and in part 
by Grant No.\ 2002068 from the United States-Israel Binational Science Foundation 
(BSF), Jerusalem, Israel}

\date{February 23, 2005} 

\begin{abstract} We prove that the Szeg\H{o} function, $D(z)$, of a measure on 
the unit circle is entire meromorphic if and only if the Verblunsky coefficients 
have an asymptotic expansion in exponentials. We relate the positions of the 
poles of $D(z)^{-1}$ to the exponential rates in the asymptotic expansion. Basically, 
either set is contained in the sets generated from the other by considering products 
of the form, $z_1 \dots z_\ell \bar z_{\ell-1}\dots \bar z_{2\ell-1}$ with $z_j$ 
in the set. The proofs use nothing more than iterated Szeg\H{o} recursion at $z$ 
and $1/\bar z$. 
\end{abstract}

\maketitle

\section{Introduction} \lb{s1} 

This paper is concerned with the spectral theory of orthogonal polynomials 
on the unit circle (OPUC) \cite{OPUC1,OPUC2,Szb,GBk1,GBk} in the case of 
particularly regular measures. Throughout, we will consider probability 
measures on $\partial\bbD=\{z\mid\abs{z}=1\}$ of the form 
\begin{equation} \lb{1.1} 
d\mu = w(\theta)\, \f{d\theta}{2\pi} + d\mu_\s 
\end{equation} 
where $w$ obeys the Szeg\H{o} condition, that is, 
\begin{equation} \lb{1.2} 
\int \log(w(\theta)) \, \f{d\theta}{2\pi} >-\infty 
\end{equation} 
In that case, the Szeg\H{o} function is defined by 
\begin{equation} \lb{1.3} 
D(z) = \exp \biggl(\int \f{e^{i\theta}+z}{e^{i\theta}-z}\, \log(w(\theta)) \, 
\f{d\theta}{4\pi}\biggr) 
\end{equation} 
Not only does $w$ determine $D$ but $D$ determines $w$, since $\lim_{r\uparrow 1} 
D(re^{i\theta})\equiv D(e^{i\theta})$ exists for a.e.~$\theta$ and 
\begin{equation} \lb{1.4} 
w(\theta) = \abs{D(e^{i\theta})}^2  
\end{equation} 
Indeed, $D$ is the unique function analytic on $\bbD =\{z\mid\abs{z}<1\}$ with $D$ 
nonvanishing on $\bbD$ so that \eqref{1.4} holds. 

Given $d\mu$, we let $\Phi_n$ be the monic orthogonal polynomial and $\varphi_n = 
\Phi_n /\|\Phi_n\|_{L^2 (d\mu)}$. The $\Phi_n$'s obey the Szeg\H{o} recursion: 
\begin{equation} \lb{1.5} 
\Phi_{n+1}(z) = z\Phi_n(z) - \bar\alpha_n \Phi_n^*(z) 
\end{equation} 
where for $P_n$ a polynomial of degree $n$, 
\begin{equation} \lb{1.6} 
P_n^*(z) = z^n \, \ol{P_n (1/\bar z)} 
\end{equation} 

The $\alpha_n$ are called {\it Verblunsky coefficients}. They lie in $\bbD$ and $\mu\mapsto 
\{\alpha_n\}_{n=0}^\infty$ is a bijection of nontrivial measures on $\partial\bbD$ 
and $\bbD^\infty$. Our goal here is to focus on the map and its inverse. Here is the 
background on our first main result: 
\begin{SL} 
\item[A.] Nevai-Totik \cite{NT89} proved that $\limsup_{n\to\infty} \abs{\alpha_n}^{1/n} 
\leq R^{-1} <1$ if and only if   
\begin{SL} 
\item[(a)] $d\mu$ obeys the Szeg\H{o} condition and $d\mu_\s =0$. 
\item[(b)] $D(z)^{-1}$ is analytic in $\{z\mid\abs{z}<R\}$. 
\end{SL} 

\item[B.] Barrios-L\'opez-Saff \cite{BLS} proved that for $R>1$, 
\begin{equation} \lb{1.7} 
\alpha_n =cR^{-n} + O(((1-\veps)R^{-1})^n) 
\end{equation} 
if and only if $D(z)^{-1}$ is meromorphic in a circle of radius $R(1+\delta)$ with a  
single, simple pole at $z=R$.  

\item[C.] Simon \cite{OPUC1} considered the functions 
\begin{equation} \lb{1.8x} 
S(z) =-\sum_{j=0}^\infty \alpha_{j-1} z^j  
\end{equation} 
(with $\alpha_{-1} =-1$) and 
\begin{equation} \lb{1.9} 
r(z) = \ol{D(1/\bar z)}\, D(z)^{-1}  
\end{equation}  
and proved that if $\limsup \abs{\alpha_n}^{1/n} \leq R^{-1} < 1$, then for some $\delta >0$, 
$r(z)-S(z)$ is analytic in $\{z\mid 1-\delta < \abs{z} <R^2 \}$ so that $S(z)$ and $r(z)$, 
which will have singularities on $\abs{z}=R$ if $\limsup \abs{\alpha_n}^{1/n} =R^{-1}$, 
must have the same singularities in $\{z\mid R\leq \abs{z} <R^2\}$. In \cite{OPUC1}, 
instead of $S(z)$ as defined by \eqref{1.8x}, one has $S(z)$ defined by $S_\book(z) = 
\sum_{j=0}^\infty \alpha_j z^j$, and the theorem is stated as analyticity of $z^{-1} 
r(z) + S_\book(z)$, equivalent to analyticity of $r-S$. But, as we will explain in 
Section~\ref{s4}, \eqref{1.8x} is the more natural object. Rather than $1-\delta 
<\abs{z} <R^2$, \cite{OPUC1} has $R^{-1} < \abs{z} < R^2$, but that is wrong since 
$\ol{D(1/\bar z)}$ can have poles at the Nevai-Totik zeros. 

\item[D.] Using Riemann-Hilbert methods, Deift-Ostensson \cite{DO} have extended the 
result on analyticity of $r(z)-S(z)$ to $\{z\mid 1-\delta <\abs{z} <R^3\}$. 

\item[E.] Barrios-L\'opez-Saff \cite{BLS2} have proven that if 
\begin{equation} \lb{1.8} 
\alpha_n =cR^{-n} + O(((1+\veps)R)^{-n-nm^2}) 
\end{equation} 
then $D(z)^{-1}$ is meromorphic in $\{z\mid\abs{z}<R^{2m-1}+\delta\}$ with poles 
precisely at $z_k =R^{2k-1}$, $k=1,2,\dots,m$. In particular, if \eqref{1.8} holds 
for all $n$, then $D(z)^{-1}$ is entire meromorphic except for poles at $R^{2k-1}$, 
$k=1,2,\dots$. 
\end{SL} 

\smallskip
Our main goal in this paper is to give a complete analysis of what can be said 
about $\alpha_n$ if $D(z)^{-1}$ is meromorphic in some disk and, contrariwise, 
about $D(z)^{-1}$ if $\alpha_n$ has an asymptotic expansion as a sum of exponentials. 
We describe our precise results in Section~\ref{s4}. 

Along the way, we found a direct, simple proof of the Deift-Ostensson result that is 
also simpler than the argument Simon used for his weaker result in \cite{OPUC1}. 
So we will give this proof next, then analyze two simple examples, and return in 
Section~\ref{s4} to a general overview and sketch of the rest of the paper. 

Of course, included among the entire meromorphic functions are the rational functions, 
and there is prior literature on this case. Szabados \cite{Sza79} considered the 
case $D(z)^{-1} = 1/q(z)$ for a polynomial $q$ and Ismail-Ruedemann \cite{IsRue} 
and Pakula \cite{Pak87} discussed $D(z)^{-1} = p(z)/q(z)$ for polynomials $p$ and 
$q$. They have some results on asymptotics of $\Phi_n$ but no discussion of links to 
the $\{\alpha_n\}_{n=1}^\infty$. As I was completing this manuscript, I received 
the latest draft of a paper of Mart\'inez-Finkelshtein, McLaughlin, and Saff \cite{MFMS} 
that has some overlap with this paper. 

\smallskip 
I would like to thank P.~Deift and J.~Ostensson for illuminating discussions.    
This research was completed during my stay as a Lady Davis Visiting Professor 
at Hebrew University, Jerusalem. I'd like to thank H.~Farkas and Y.~Last for 
the hospitality of the Mathematics Institute at Hebrew University.

\section{The $R^3$ Result} \lb{s2} 

Our goal in this section is to prove 

\begin{theorem}\lb{T2.1} Let 
\begin{equation} \lb{2.1} 
\limsup \abs{\alpha_n}^{1/n} = R^{-1} < 1 
\end{equation} 
so that $D(z)^{-1}$ and $S(z)$ are analytic in $\{z\mid\abs{z}<R\}$. Then for some 
$\delta >0$, $r(z) -S(z)$ is analytic in $\{z\mid 1-\delta <\abs{z} <R^3\}$. 
\end{theorem} 

\begin{remarks} 1. Here $D$ is given by \eqref{1.3}, $S(z)$ by \eqref{1.8x}, and 
$r(z)$ by \eqref{1.9}. 

2. The proof will make repeated use of Szeg\H{o} recursion \eqref{1.5}. 
\end{remarks} 

We introduce the symbol $\wti O$ by $f=\wti O(g)$ if and only if for all $\veps$, 
$\abs{f}/\abs{g}^{1-\veps} \to 0$. 

\begin{lemma}\lb{L2.2} Let \eqref{2.1} hold. Then 
\begin{SL} 
\item[{\rm{(a)}}] For all $\veps >0$, 
\begin{equation} \lb{2.2} 
\sup_{n,\, \abs{z}<R-\veps}\, \abs{\Phi_n^*(z)} <\infty  
\end{equation} 

\item[{\rm{(b)}}] For $\abs{z}\leq 1$, 
\begin{equation} \lb{2.3} 
\abs{\Phi_n(z)} =\wti O(\max (R^{-1}, \abs{z})^n) 
\end{equation} 

\item[{\rm{(c)}}] For $\abs{z}\leq 1$, 
\begin{equation} \lb{2.4} 
\abs{\Phi_n^*(z) -D(0)D(z)^{-1}} =\wti O(R^{-n} \max (R^{-1}, \abs{z})^n) 
\end{equation} 
\end{SL} 
\end{lemma} 

\begin{remarks} 1. There is an implicit uniformity in $z$ in the $\wti O$ statements 
\eqref{2.3}, \eqref{2.4}. 

2. (a) is due to Nevai-Totik; (b) appears in Simon \cite{Saff1,Saff2}. 
\end{remarks}

\begin{proof} (a) \ From \eqref{1.5} and $\abs{\Phi_n (e^{i\theta})} = 
\abs{\Phi_n^* (e^{i\theta})}$, we see 
\begin{equation} \lb{2.5} 
\sup_{\abs{z}=1}\, \abs{\Phi_n(z)} \leq \prod_{j=0}^{n-1} (1+\abs{\alpha_j})  
\end{equation} 
so, by $\prod_{j=0}^\infty (1+\abs{\alpha_j})<\infty$ and the maximum principle, 
\begin{equation} \lb{2.6} 
\sup_{n,\, \abs{z}\leq 1}\,  \abs{\Phi_n^*(z)} <\infty 
\end{equation} 
from which we get, by \eqref{1.6}, that 
\begin{equation} \lb{2.7} 
C_1 \equiv \sup_{n,\, \abs{z}\geq 1}\, \abs{z}^{-n} \abs{\Phi_n(z)} <\infty  
\end{equation} 

The $\,^*\,$ of \eqref{1.5} is 
\begin{equation} \lb{2.8} 
\Phi_{n+1}^*(z) - \Phi_n^*(z) =-\alpha_n z\Phi_n(z)  
\end{equation} 
so that 
\begin{equation} \lb{2.9} 
\Phi_n^*(z) = 1-\sum_{j=0}^{n-1} \alpha_j z \Phi_j(z)  
\end{equation} 
Thus, by \eqref{2.7}, 
\begin{equation} \lb{2.10}  
\abs{\Phi_n^*(z)} \leq  1 + C_1 \sum_{j=0}^{n-1}\, \abs{\alpha_j}\, \abs{z}^{j+1}  
\end{equation} 
Given \eqref{2.1}, we see that \eqref{2.2} holds. 

\smallskip 
(b) \ \eqref{2.2} and \eqref{1.6} imply that for $\abs{z}>R^{-1} + \veps$, 
\begin{equation} \lb{2.11} 
\abs{\Phi_n(z)} \leq C_\veps \abs{z}^n \qquad (\abs{z} >R^{-1} +\veps)  
\end{equation}  
This plus the maximum principle implies \eqref{2.3}. 

\smallskip 
(c) \ It is a theorem of Szeg\H{o} \cite{Sz20,Sz21} (see Theorem~2.4.1 of \cite{OPUC1}) 
that in $\abs{z} <1$, 
\begin{equation} \lb{2.12} 
\lim_{n\to\infty}\, \Phi_n^*(z) =D(0) D(z)^{-1} \equiv d(z)^{-1}   
\end{equation} 
Thus, summing \eqref{2.8} to infinity, 
\begin{equation} \lb{2.13} 
\abs{d(z)^{-1} -\Phi_n^*(z)} \leq \sum_{j=n}^\infty\, \abs{\alpha_j} \, \abs{z}\, 
\abs{\Phi_j (z)}  
\end{equation} 
Since $\alpha_j =\wti O(R^{-j})$ and \eqref{2.3} holds, we obtain \eqref{2.4}. 
\end{proof} 

\begin{proof}[Proof of Theorem~\ref{T2.1}] We use the function $d(z)$ of \eqref{2.12}. 
Since $\Phi_n^*(z) \to d(z)^{-1}$ for $\abs{z}<1$ and \eqref{1.2} holds, 
the Vitali theorem implies $d(z)^{-1}$ is analytic in $\{z\mid\abs{z}<R\}$ and 
$\Phi_n^*(z)\to d(z)^{-1}$ in that region. By summing \eqref{2.8} to infinity, 
\begin{equation} \lb{2.14} 
d(z)^{-1} =1-\sum_{j=1}^\infty \alpha_{j-1} z\Phi_{j-1} (z) 
\end{equation} 
which we write 
\begin{equation} \lb{2.15} 
d(z)^{-1} = \ol{d(1/\bar z)}^{-1} S(z) + \bigl[ 1-\ol{d(1/\bar z)}^{-1}\bigr] - 
\sum_{j=1}^\infty \alpha_{j-1} z \bigl[ \Phi_{j-1}(z) - \ol{d(1/\bar z)}^{-1}\, 
z^{j-1}\bigr]  
\end{equation} 
where this formula is valid in $\{z\mid R^{-1} <\abs{z} <R\}$. 

Apply $\,^*\,$ to \eqref{2.4} and see that in $\abs{z}\geq 1$, 
\begin{equation} \lb{2.16} 
\abs{\Phi_n(z) -\ol{d(1/\bar z)}^{-1} z^n} \leq \wti O(R^{-n} \max (\abs{z}R^{-1}, 1)^n) 
\end{equation} 
Thus the summand in \eqref{2.15} is bounded in $\{z\mid\abs{z}\geq 1\}$ by $\abs{z} 
\wti O(R^{-2n} \max (\abs{z}R^{-1}, 1)^n)$. In $1\leq\abs{z}\leq R$, this is bounded by 
$\wti O(RR^{-2n})$ and in $R \leq \abs{z}$ by $\wti O(\abs{z}^{n+1} R^{-3n})$. Thus, the 
sum in \eqref{2.15}, which is a sum of functions each analytic in $\{z\mid\abs{z} >R^{-1}\}$, 
converges uniformly in $\{z\mid 1 \leq \abs{z} <R^3\}$. Multiplying by $\ol{d(1/\bar z)}$,  
which is analytic in $\{z\mid\abs{z} > 1-\delta\}$, implies the result. 
\end{proof}

\section{Two Examples} \lb{s3}  

We want to analyze two examples from \cite{OPUC1} from the point of view of singularities 
of $D(z)^{-1}$ and asymptotics of $\alpha_n$. The first is already mentioned in this 
context in \cite{BLS2}. 

\begin{example}[Rogers-Szeg\H{o} polynomials; Example~1.6.5 of \cite{OPUC1}] \lb{E3.1} 
Here $0<q<1$, 
\begin{equation} \lb{3.1} 
\alpha_n = (-1)^n q^{(n+1)/2} 
\end{equation} 
and 
\begin{equation} \lb{3.2} 
D(z) =\prod_{j=0}^\infty (1-q^j)^{1/2} (1+q^{j+1/2} z) 
\end{equation}
Let $R=q^{-1/2}$. Then 
\begin{align} 
S(z) &= -\sum_{j=0}^\infty (-1)^{j-1} q^{j/2} z^j \notag \\
&= (1+zR^{-1})^{-1} \lb{3.3} 
\end{align} 
has a single pole at 
\begin{equation} \lb{3.4} 
z_1 = -R 
\end{equation} 
On the other hand, by \eqref{3.2}, $D(z)$ has a zero and so $D(z)^{-1}$ a pole at 
\begin{equation} \lb{3.5} 
z_\ell = -R^{2\ell -1} \qquad \ell=1,2,\dots  
\qedhere
\end{equation} 
\end{example}  

\begin{example}[Single nontrivial moment; Example~1.6.4 of \cite{OPUC1}] \lb{E3.2}
Fix $0<a<1$ and let 
\begin{equation} \lb{3.6} 
d\mu_a (\theta) = (1-a\cos\theta)\, \f{d\theta}{2\pi} 
\end{equation} 
Let 
\begin{equation} \lb{3.7} 
\mu_\pm = \f{1}{a} \pm \sqrt{\biggl(\f{1}{a}\biggr)^2 -1} 
\end{equation} 
so $\mu_-\mu_+ =1$ and $\mu_- <1$. Then 
\begin{align} 
D(z) &= \sqrt{\f{a}{2\mu_-}}\, (1-\mu_- z) \lb{3.8} \\ 
&= \sqrt{\f{a}{2\mu_-}}\, \biggl(1 -\f{z}{\mu_+}\biggr) \lb{3.9} 
\end{align}
so $D(z)^{-1}$ has a single pole at 
\begin{equation} \lb{3.10} 
z_1 = \mu_+ 
\end{equation}
On the other hand, 
\begin{align} 
\alpha_n &= \f{-(\mu_+ - \mu_-)}{(\mu_+^{n+2} - \mu_-^{(n+2)})} \notag \\
&= -(\mu_+ -\mu_-)\mu_+^{-n-2} (1-\mu_+^{-(2n+4)})^{-1} \notag \\
&= -(\mu_+ -\mu_-) \sum_{j=1}^\infty (\mu_+^{-n-2})^{2j-1} \lb{3.11}
\end{align} 
and so $S(z)$ has  poles at 
\begin{equation} \lb{3.12} 
z_j = \mu_+^{2j-1} \qquad j=1,2,\dots 
\end{equation}
\end{example}  

In these examples, the set of singularities of $S$ and of $D(z)^{-1}$ are distinct 
and one or the other might be larger. If $\{z_j\}$ are the singularities, then 
$\{z_j^{2\ell -1}\}_{j,\ell=1,2,\dots}$ are identical for $S$ and $D^{-1}$, which 
motivates the $\bbG$ construction of the next section.

\section{Overview and Discussion of Further Results} \lb{s4} 

\begin{definition} A sequence $\{A_n\}_{n=-1}^\infty$ of complex numbers is said to have 
an asymptotic series with error $R^{-n}$ for some $R>1$ if and only if there exists a 
finite number of points $\{\mu_j\}_{j=1}^J$ in $\{w\mid 1<\abs{w}<R\}$ and polynomials 
$\{P_j\}_{j=1}^J$ so that 
\begin{equation} \lb{4.1} 
\limsup_{n\to\infty} \, \biggl| A_n -\sum_{j=1}^J P_j(n) \mu_j^{-(n+1)}\biggr|^{1/n} 
\leq R^{-1} 
\end{equation} 
Equivalently, 
\[
A_n =\sum_{j=1}^J P_j(n) \mu_j^{-n-1} + \wti O(R^{-n}) 
\]
We say $A_n$ has a {\it complete asymptotic series\/} if it has an asymptotic series 
with error $R^{-n}$ for all $R>1$. 
\end{definition} 

In many ways, our main result in this paper is: 

\begin{theorem} \lb{T4.1} Let $d\mu$ be a nontrivial probability measure on $\partial\bbD$ 
with Verblunsky coefficients, $\alpha_n$. Then $\alpha_n$ has a complete asymptotic series 
if and only if 
\begin{SL} 
\item[{\rm{(1)}}] $d\mu_\s =0$ and $d\mu$ obeys the Szeg\H{o} condition. 
\item[{\rm{(2)}}] $D(z)^{-1}$ is an entire meromorphic function. 
\end{SL} 
\end{theorem} 

Of course, 
\begin{equation} \lb{4.2} 
\sum_{n=0}^\infty z^n \mu_j^{-n} = \biggl( 1-\f{z}{\mu_j}\biggr)^{-1}  
\end{equation} 
and so, taking derivatives, for $\ell=1,2\dots $, 
\begin{equation} \lb{4.3} 
\sum_{n=0}^\infty (n+\ell)(n+\ell-1) \dots (n+1)z^n \mu_j^{-n} 
= \ell! \bigg(1 - \f{z}{\mu_j}\biggr)^{-\ell-1} 
\end{equation} 
So \eqref{4.1} is equivalent to a sum of explicit pole terms:  

\begin{proposition} \lb{P4.2} $\{A_n\}_{n=-1}^\infty$ has an asymptotic series with error 
$R^{-n}$ if and only if 
\begin{equation} \lb{4.4} 
F(z) =\sum_{n=0}^\infty A_{n-1} z^n 
\end{equation} 
is meromorphic in $\{z\mid \abs{z}<R\}$ with a finite number of poles, all in $\{z\mid 
1< \abs{z}<R\}$. In particular, $\{A_n\}_{n=-1}^\infty$ has a complete asymptotic series 
if and only if $F(z)$ is an entire meromorphic function. 
\end{proposition} 

Thus, Theorem~\ref{T4.1} is equivalent to 

\begin{theorem}\lb{T4.3} {\rm{(1)}} and {\rm{(2)}} of Theorem~\ref{T4.1} are equivalent 
to the function $S$ of \eqref{1.8x} being an entire meromorphic function. 
\end{theorem} 

Note that the $\mu_j$'s and $P_j$'s are determined uniquely by the $A_n$'s. 

\smallskip
Both to prove the results and for its intrinsic interest, we are interested in the 
relation between the poles of $S(z)$ and of $D(z)^{-1}$ and in results in fixed 
circles. By a {\it discrete exterior set}, we mean a subset, $T$, of $\{w\mid 
1<\abs{w} <\infty\}$ so that $\#[\{w\mid 1<\abs{w} <R\}\cap T]$ is finite for 
each $R>1$. Given a discrete exterior set $T$, define for $k=1,2,\dots$,  
\begin{align} 
\bbG^{(2k-1)} (T) &= \{\lambda_{i_1} \dots \lambda_{i_k} 
\bar\lambda_{i_{k+1}} \dots \bar\lambda_{i_{2k-1}} \mid \lambda_j \in T\}\lb{4.5} \\ 
\bbG_{2k-1}(T) &= \bigcup_{j=k}^\infty \, \bbG^{(2j-1)}(T) \lb{4.6} \\ 
\bbG(T) &= \bbG_1(T) \lb{4.7} 
\end{align} 
$\bbG(T)$ will be called the generated set. Note that 
\begin{equation} \lb{4.8} 
\bbG(\bbG(T)) = \bbG(T) 
\end{equation} 

We will prove the following: 

\begin{theorem} \lb{T4.4} Let $\alpha_n$ be a set of Verblunsky coefficients with 
complete asymptotic series and let $T$ be the set of $\lambda_j$'s that enter in 
the series. Let $P$ be a set of poles of $D(z)^{-1}$. Then 
\begin{equation} \lb{4.9} 
T\subset \bbG(P) \quad\text{and}\quad P\subset \bbG(T) 
\end{equation}
\end{theorem} 

This implies the following refined form of Theorem~\ref{T4.1}: 

\begin{theorem}\lb{T4.5} Let $Q$ be an exterior discrete set with $\bbG(Q) =Q$. Then 
$\alpha_n$ is a set of Verblunsky coefficients with $\lambda_j$'s in $Q$ if and only 
if condition {\rm{(1)}} of Theorem~\ref{T4.1} holds and the poles of $D(z)^{-1}$ lie 
in $Q$. 
\end{theorem}  

Theorems~\ref{T4.1}, \ref{T4.4}, and \ref{T4.5} are equivalence results and thus 
both a direct (going from $\alpha$ to $D$) and inverse (going from $D$ to $\alpha$) 
aspect. Generally, direct arguments are simpler than inverse. We will actually deduce 
everything from direct arguments and a bootstrap. An inverse argument is only used 
to start the analysis, and that was already done by Nevai-Totik. Here is the 
master stepping stone we will need. Throughout, we suppose there is $R>1$ so 
\begin{equation} \lb{4.10} 
\limsup_{n\to\infty}\, \abs{\alpha_n}^{1/n} = R^{-1}  
\end{equation} 

\begin{theorem}\lb{T4.6} Fix $\ell=1,2,\dots$. Suppose $S(z)$ is meromorphic in 
\begin{equation} \lb{4.11x} 
\calR_\ell = \{z\mid 0<\abs{z} <R^{2\ell-1}\} 
\end{equation} 
Then $D(z)^{-1}$ is meromorphic there and the poles of $D(z)^{-1}$ there lie in 
$\bbG(T_\ell)$ where $T_\ell$ is the set of poles of $S(z)$ in $\calR_\ell$. 
Moreover, $r(z)-S(z)$ has a meromorphic continuation to $\calR_{\ell+1} \cap 
\{z\mid 1 > \abs{z} -\delta\}$ and the poles of this difference lie in $\bbG_3 (T_\ell)$.  
\end{theorem} 

We are heading towards a proof that Theorem~\ref{T4.6} implies the earlier 
Theorems~\ref{T4.1}, \ref{T4.4}, and \ref{T4.5}. We need a preliminary notion and 
fact. 

\begin{definition} Let $Q$ be an exterior discrete set with $\bbG(Q) =Q$. We say 
that $W\subset Q$ is a set of minimal generators if and only if $\bbG(W)=Q$ and 
$\bbG_3(W)\cap W=\emptyset$. 
\end{definition} 

\begin{proposition}\lb{P4.7} Any exterior discrete set $Q$ with $\bbG(Q)=Q$ has a 
minimal set of generators. 
\end{proposition} 

\begin{proof} Order the points in $Q$, $w_1, w_2, \dots$ so $\abs{w_n} \leq  
\abs{w_{n+1}}$. Define $W$ inductively by putting $w_n$ in $W$ if and only if 
$w_n\notin \bbG_3 (\{w_1, \dots, w_{n-1}\})$. It is easy to see that $W$ is a set 
of minimal generators.  
\end{proof} 

\begin{proof}[Proof that Theorem~\ref{T4.6} implies Theorems~\ref{T4.4}, 
\ref{T4.5}, and \ref{T4.1}] It suffices to prove that $D^{-1}$ is entire meromorphic 
if and only if $S$ is, and to prove Theorem~\ref{T4.4} since it in turn implies 
Theorem~\ref{T4.5}, which implies Theorem~\ref{T4.1}. If $S(z)$ is entire meromorphic, 
it is meromorphic in each $\calR_\ell$, so $D(z)^{-1}$ is meromorphic in each 
$\calR_\ell$ and, clearly, $P\subset \bbG(T)$. 

Conversely, if $D(z)^{-1}$ is entire meromorphic, we prove $S(z)$ is entire 
meromorphic by proving inductively that it is meromorphic in each $\calR_\ell$. 
$S(z)$ is meromorphic in $\calR_1$ by the Nevai-Totik theorem. If we know $S(z)$ 
is meromorphic in $\calR_\ell$, then by Theorem~\ref{T4.6}, $r(z)-S(z)$ is 
meromorphic in $\calR_{\ell+1}\backslash\calR_\ell$ so, since $r(z)$ is meromorphic 
on $\calR_{\ell+1}$, we conclude that $S(z)$ is meromorphic there also. 

Finally, to identify the points of $T$, as lying in $\bbG(P)$ with $P$ the poles of 
$D^{-1}(z)$, suppose $W$ is a set of minimal generators of $T$. If $w_j\in W$\!, then 
$w_j\notin \bbG_3 (T)$, so $S-r$ is regular at $w_j$ by Theorem~\ref{T4.6}. Since 
$w_j$ is a singularity of $S$, it must be a singularity of $r$, that is, $w_j \in P$. 
Thus, $T=\bbG(W)\subset \bbG(P)$. 
\end{proof} 

Our proof of Theorem~\ref{T4.6} will also show 

\begin{theorem}\lb{T4.8} Suppose $z_0\in\bbG^3 (T)$ has a unique expression as 
$z_0 = \mu_1^2 \bar\mu_2$ with $\mu_1,\mu_2\in T$. Suppose also that $z_0\notin 
T\cup\bbG_5 (T)$. Then $r(z)-S(z)$ has a singularity at $z_0$. 
\end{theorem}  

\begin{remarks} 1. For example, if $S(z)$ has a single pole, $z_0$, with $\abs{z_0} 
=R$, then either $S$ or $D^{-1}$ or both have a pole at $z_0 R^2$. 

2. Our proof shows that if the poles of $S(z)$ have $\{\log \abs{z_j}\}$ independent 
over the rationals, $D^{-1}$ has a pole at every point in $\bbG(T)$. 

3. Our proof also allows the precise calculation of the singularity in $S-r$ at 
any point in $\bbG(T)$. There can be cancellations if $z_0$ can be written as a 
product in $\bbG(T)$ in more than one way. So one cannot guarantee a singularity 
of $r(z)-S(z)$ at every point of $\bbG_3(T)$, but that will happen in some 
generic sense. 
\end{remarks}

Note that our results generalize those of Barrios-L\'opez-Saff \cite{BLS2} in three ways: 
\begin{SL} 
\item[(a)] They only have results on the the direct problem, that is, going from 
$\alpha$ to $D^{-1}$, while we have results in both directions. 
\item[(b)] They allow only a single term in the $\alpha$ asymptotics. 
\item[(c)] Their error assumptions in case of disks are much stronger ($R^{-nm^2}$ 
vs.\ $R^{-n(2m+1)}$) than ours. 
\end{SL} 

This concludes the description of our main results --- and reduces everything to 
proving Theorems~\ref{T4.6} and \ref{T4.8}. We will do this for $2\ell -1=3$ 
in Section~\ref{s5} and general $2\ell-1$ in Section~\ref{s6}. 

One could analyze other situations such as where $S(z)$ has a branch cut associated 
with specific asymptotics for $\alpha_n$ such as $n^\beta R^{-n}$ with $\beta$ 
nonintegral. 

We close this section, which is the continuation of the introduction, with two remarks. 
First, there is a scattering theoretic interpretation of $S$ and $r$. Since 
$\sum_{n=0}^\infty \abs{\alpha_n}^2 <\infty$, one can define wave operators  
(see Geronimo-Case \cite{GC} and Section~10.7 of \cite{OPUC2}), $\Omega^\pm : 
L^2 (\partial\bbD, \f{d\theta}{2\pi})\to L^2 (\partial\bbD, \abs{D}^2 
\f{d\theta}{2\pi})$, which obey  
\begin{equation} \lb{4.11} 
(\Omega^+ f) (\theta) = D(e^{i\theta})^{-1} f(\theta) \qquad 
(\Omega^- f)(\theta) = \ol{D(e^{i\theta})}^{-1} f(\theta) 
\end{equation} 
Thus, the reflection coefficient is given by 
\begin{equation} \lb{4.12} 
((\Omega^-)^{-1} \Omega^+f)(\theta) = \ol{D(e^{i\theta})}\, D(e^{i\theta})^{-1} 
f(\theta)
\end{equation} 
so $r(z)$ is the analytic continuation of the reflection coefficient. $S(z)$ is 
the leading Born approximation to $r$ (see Newton \cite{Newt} and Chadan-Sabatier 
\cite{ChaSab} for background on scattering theory). While we will not study it 
from this point of view, it is presumably true that the arguments in the next 
two sections can be interpreted as use of some kind of Born series. 

The second issue concerns a comparison between the basic formula used by 
Nevai-Totik \cite{NT89} to do the inverse problem and a different, but 
similar-looking, formula used in our discussion, namely \eqref{2.14}. The formula 
they use, where they quote Freud \cite{FrB}, is also in Geronimus \cite{GBk}: 
\begin{equation} \lb{4.13} 
\alpha_n = -\kappa_\infty \int \ol{\Phi_{n+1} (e^{i\theta})}\, D(e^{i\theta})^{-1} \, 
d\mu(\theta)  
\end{equation} 
where $\kappa_\infty =\lim_{n\to\infty}\kappa_n$ with $\kappa_n =\|\Phi_n\|^{-1}$,  
so $\kappa_n = \varphi_n^*(0)$ and 
\begin{equation} \lb{4.14} 
\kappa_\infty = D(0)^{-1}   
\end{equation} 
\eqref{4.13} only holds if $d\mu_\sing =0$. 

Since $\varphi_n =\kappa_n \Phi_n$, \eqref{4.13} can be rewritten as 
\begin{equation} \lb{4.15} 
\alpha_n = -\kappa_\infty \kappa_n^{-1} \langle \varphi_{n+1}, D^{-1} \rangle  
\end{equation} 
Since 
\[ 
\langle 1,D^{-1}\rangle =\int \ol D(e^{i\theta}) \, \f{d\theta}{2\pi} = 
\ol{D(0)} = \kappa_\infty^{-1} 
\]
\eqref{4.15} also holds if we interpret $\alpha_{-1}=-1$ and $\kappa_{-1}=1$. 
Thus, with $\alpha_{-1}=-1$, \eqref{4.15} is equivalent to 
\begin{equation} \lb{4.16} 
d(z)^{-1} \equiv D(0) D(z)^{-1} = - \kappa_\infty^{-2} \sum_{n=-1}^\infty 
\kappa_n \alpha_n \varphi_n(z)  
\end{equation} 

On the other hand, \eqref{2.14} says 
\begin{equation} \lb{4.17} 
(d(z)^{-1}-1) z^{-1} =-\sum_{n=0}^\infty \alpha_n \kappa_n^{-1} \varphi_n(z)
\end{equation} 
or equivalently, 
\begin{equation} \lb{4.18} 
\alpha_n =-\kappa_\infty^{-1} \kappa_n^2 \int \ol{\Phi_n (e^{i\theta})}\, 
[D(e^{i\theta})^{-1} - D(0)^{-1}] e^{-i\theta}\, d\mu(\theta)  
\end{equation} 
These formulae are distinct, and it is striking that both are true and their proofs 
(see (2.4.35) of \cite{OPUC1}) are so different. Where Nevai-Totik \cite{NT89} use 
\eqref{4.13}, one could just as well use \eqref{4.18}.

\section{The $R^5$ Result} \lb{s5} 

In this section, as a warmup and also as the start of induction for the general case, we 
consider the case $2\ell -1=3$, that is, $\ell=2$ where we deal with induced  
singularities in $\{z\mid R^3 \leq \abs{z} < R^5\}$. Thus, we should suppose 
\begin{equation} \lb{5.1} 
\alpha_n = \sum_{k=1}^K P_k (n) \mu_k^{-n-1} + \wti O(R^{-3n})  
\end{equation} 
with $R\leq \abs{\mu_k} < R^3$. Here, $P_k(n)$ are polynomials. We will instead 
suppose that 
\begin{equation} \lb{5.2} 
\alpha_n = \sum_{k=1}^K c_k \mu_k^{-n-1} + \wti O(R^{-3n}) 
\end{equation} 
The consideration of general $P_k$'s rather than constants presents no difficulties 
other than notational ones, so we spare the reader. Our goal is to prove 

\begin{theorem}\lb{T5.1} If \eqref{5.2} holds, then $D(z)^{-1}$ is meromorphic in 
$\{z\mid\abs{z} <R^3\}$ with poles precisely at $\{\mu_k\}_{k=1}^K$. In addition, $S(z) 
-r(z)$ is meromorphic in $\{z\mid 1-\delta < \abs{z} < R^5\}$ with poles contained 
in $\bbG^{(3)}(\{\mu_k\}_{k=1}^K)$. Moreover, if $z_0 = \mu_{i_1}^2 \bar\mu_{i_2}$ 
in precisely one way and $\abs{z_0} \leq R^5$, then $S(z)-r(z)$ has a pole at $z_0$. 
\end{theorem} 

We note that the first statement is immediate from Theorem~\ref{T2.1}, so we will 
focus on $R^3 \leq \abs{z} < R^5$. We will follow the same three-step strategy used 
in Section~\ref{s2}: 
\begin{SL} 
\item[(i)] Estimate $\Phi_n$ in $\{z\mid\abs{z} < R^{-3} (1+\delta)\}$. 
\item[(ii)] Estimate $\Phi_n^* -d(z)^{-1}$ in $\{z\mid\abs{z} < R^{-3} (1+\delta)\}$ 
using \eqref{2.8} and Step (i). 
\item[(iii)] Estimate $S(z) -r(z)$ in $\abs{z} > R^3/(1+\delta)$ using \eqref{2.8}, 
the formula $\Phi_n(z) z^{-n} = \ol{\Phi_n^*(1/\bar z)}$, and the estimate in 
Step (ii).  
\end{SL} 

\smallskip 
What will be different from Section~\ref{s2} is that we will find the leading 
asymptotics of $\Phi_n$ rather than just use $\abs{\Phi_n(z)} \leq \wti O(R^{-n})$ 
in $\abs{z}<R$. In essence, this leading asymptotics was discussed in \cite{Saff2}  
and we will use the techniques from there, although in a slightly more general 
context. 

\begin{theorem} \lb{T5.2} Suppose \eqref{5.2} holds. Choose $\delta$, so for all $k$, 
$\abs{\mu_k} <R^3 /(1+\delta)$. Define $E_n(z)$ in $\{z\mid\abs{z} < R^{-3}(1+\delta)\}$ 
by 
\begin{equation} \lb{5.3} 
\Phi_n(z) = - d(z)^{-1} \sum_{k=1}^K \bar c_k \bar\mu_k^{-n} (1- z\bar\mu_k)^{-1} + E_n(z) 
\end{equation} 
Then, for $\abs{z} < R^{-3} (1+\delta)$, 
\begin{equation} \lb{5.4} 
\abs{E_n(z)} = \wti O(\max (R^{-3}, \abs{z})^n) 
\end{equation} 
\end{theorem} 

\begin{remark} Since $\abs{\mu_k} <R^3 (1+\delta)^{-1}$ and $\abs{z}<R^{-3}(1+\delta)$, 
$\abs{\abs{z} \mu_k}<1$. 
\end{remark} 

\begin{proof} Iterating \eqref{1.5} from $j=n-1$ down to $j=0$ yields 
\begin{equation} \lb{5.5} 
\Phi_n(z) =z^n - \sum_{j=1}^n \bar\alpha_{n-j} z^{j-1} \Phi_{n-j}^*(z)  
\end{equation} 
Write \eqref{5.2} as 
\begin{align} 
\alpha_n &= \sum_{k=1}^K c_k \mu_k^{-n-1} + (\delta\alpha)_n \lb{5.6} \\
(\delta\alpha)_n &= \wti O(R^{-3n}) \lb{5.7}  
\end{align} 

In \eqref{5.5}, do the following: 
\begin{align} 
E_0^{(n)} &= z^n \lb{5.8} \\
E_1^{(n)} &= -\sum_{j=1}^n \bar\alpha_{n-j} z^{j-1} [\Phi_{n-j}^*(z) -d(z)^{-1}] \lb{5.9} \\ 
E_2^{(n)} &= -\sum_{j=1}^n \,\ol{(\delta\alpha)}_{n-j}\, z^{j-1}\, d(z)^{-1} \lb{5.10} \\
E_3^{(n)} &= \sum_{k=1}^K \bar c_k \sum_{j=n+1}^\infty \bar\mu_k^{-(n+1-j)} z^{j-1}\, 
d(z)^{-1} \lb{5.11} 
\end{align} 
Since $\sum_{j=1}^\infty \bar\mu_k^{-(n+1-j)} z^{j-1} = \bar\mu_k^{-1} (1-z\bar\mu_k)^{-1}$, 
\eqref{5.3} holds where 
\[
E_n(z) = E_0^{(n)} + E_1^{(n)} + E_2^{(n)} + E_3^{(n)} 
\] 

We need to show that for $j=0,1,2,3$, 
\begin{equation} \lb{5.12} 
\abs{E_j^{(n)} (z)} =\wti O(\max(R^{-3},\abs{z})^n) 
\end{equation} 
This is trivial for $j=0$. For $j=1$, we use \eqref{2.4} and $\abs{\alpha_{n-j}} 
=\wti O(R^{-(n-j)})$ to see if $\abs{z}<R^{-1}$:  
\begin{align*} 
\abs{E_1^{(n)}} &\leq \sum_{j=1}^n \, \abs{z}^{j-1} \wti O(R^{-(n-j)}) \wti O(R^{-2(n-j)}) \\
&= n\wti O( \max (\abs{z}, R^{-3})^n) \\ 
&= \wti O(\max(\abs{z}, R^{-3})^n) 
\end{align*} 

By \eqref{5.7}, $d(z)^{-1}$ is bounded in $\{z\mid\abs{z}<R^{-1}\}$ and  
$n\wti O(\max(\abs{z}, R^{-3})^n) =\wti O(\max (\abs{z}, R^{-3})^n)$, we have that
\eqref{5.12} holds for $j=2$. For $j=3$, we note that the sum of the geometric 
series is $z^n (1-\bar\mu_k z)^{-1}$, so in $\abs{z} < R^{-3} (1+\delta)$, 
we get a $\abs{z}^n$ bound. This proves \eqref{5.12}. 
\end{proof} 

\begin{theorem}\lb{T5.3} Let \eqref{5.2} hold and let $\delta$ be as in Theorem~\ref{T5.2}. 
Define $\wti E_n$ for $\abs{z} <R^{-3} (1+\delta)$,  
\begin{equation} \lb{5.13} 
\Phi_n^*(z) -d(z)^{-1} = -\biggl[d(z)^{-1} \sum_{k,\ell =1}^K \bar c_k c_\ell \mu_\ell^{-1} 
(1-z\bar\mu_k)^{-1} (1-\bar\mu_k^{-1} \mu_\ell^{-1}) (\bar\mu_k \mu_\ell)^{-n}\biggr] 
+\wti E_n(z)  
\end{equation} 
Then for $\abs{z} <R^{-3} (1+\delta)$, 
\begin{equation} \lb{5.14} 
\wti E_n(z) = \wti O(R^{-n} \max (\abs{z}, R^{-3})^n)  
\end{equation} 
\end{theorem} 

\begin{proof} We iterate \eqref{2.8} to get 
\begin{equation} \lb{5.15} 
\Phi_n^*(z) - d(z)^{-1} = \sum_{j=n}^\infty \alpha_j z\Phi_j(z) 
\end{equation} 
In \eqref{5.15}, first replace $\alpha_j$ by \eqref{5.6} and then, in the main term, replace 
$\Phi_j$ by \eqref{5.3}. Noting that $\sum_{j=n}^\infty \bar\mu_k^{-j} \mu_\ell^{-j-1} = 
\mu_\ell^{-1} (1-\bar\mu_k^{-1}\mu_\ell^{-1})^{-1} (\bar\mu_k\mu_\ell)^{-n}$, we see that 
\eqref{5.13} holds, where 
\[
\wti E_n = \wti E_1^{(n)} + \wti E_2^{(n)} 
\]
with 
\begin{align} 
\wti E_1^{(n)} &= \sum_{j=n}^\infty \, (\delta\alpha_j) z\Phi_j(z) \lb{5.16} \\ 
\wti E_2^{(n)} &= \sum_{j=n}^\infty  \biggl( \, \sum_{k=1}^K c_k \mu_k^{-j-1}\biggr) 
zE_j(z) \lb{5.17} 
\end{align}

By \eqref{5.7} and \eqref{2.3}, 
\begin{equation} \lb{5.18} 
\abs{E_1^{(n)}} \leq c \sum_{j=n}^\infty \wti O(R^{-3n}) \wti O(R^{-n}) 
= \wti O(R^{-4n}) 
\end{equation} 
By \eqref{5.4}, 
\begin{align} 
\abs{E_2^{(n)}} &\leq c \sum_{j=n}^\infty R^{-j-1}\, \wti O(\max(R^{-3}, \abs{z})^j) \lb{5.19} \\ 
&= \wti O(R^{-n} \max(\abs{z},R^{-3})^n) \notag
\end{align} 
proving \eqref{5.14}. 
\end{proof} 

\begin{theorem}\lb{T5.4} Suppose \eqref{5.2} holds and $\delta$ is chosen as in 
Theorem~\ref{T5.2}. Then in $\{z\mid R^3 (1+\delta)^{-1} <\abs{z} <R^5\}$, 
we have that $r(z)-S(z) -q_3(z)$ is analytic where 
\begin{equation} \lb{5.20} 
q_3(z) = - \sum_{k,\ell,r=1}^K c_k\bar c_\ell c_r\, 
z(z-\mu_k)^{-1} (1-\mu_k^{-1} \bar\mu_\ell^{-1})^{-1} \mu_k (1-z\mu_k^{-1} 
\bar\mu_\ell^{-1} \mu_r^{-1})^{-1} 
\end{equation} 
\end{theorem} 

\begin{proof} By \eqref{2.15}, 
\begin{equation} \lb{5.21} 
r(z)-S(z) = \bigl[\,\ol{d(1/\bar z)}-1] - \ol{d(1/\bar z)}\, \sum_{j=1}^\infty 
\alpha_{j-1} z[\Phi_{j-1} -\ol{d(1/\bar z)}^{-1} z^{-j-1}\bigr] 
\end{equation} 
Because $q_3(z)$ is obtained by summing 
\begin{equation} \lb{5.22} 
\sum_{j=1}^\infty \, (\alpha_{j-1} -\delta\alpha_{j-1}) \bigl[ \Phi_{j-1} - 
\ol{d(1/\bar z)}\, z^{j-1} - z^{j-1}\, \ol{\wti E_{j-1} (1/\bar z)}\,\bigr]
\end{equation} 
we see that 
\begin{equation} \lb{5.23} 
r(z)-S(z) -q_3(z) = -\ol{d(1/\bar z)}\, \sum_{j=1}^\infty F_{1,j}(z) + F_{2,j}(z)  
\end{equation} 
where 
\begin{align} 
F_{1,j}(z) &= \delta\alpha_{j-1} z [\Phi_{j-1} -d(1/\bar z) z^{j-1}] \lb{5.24} \\ 
F_{2,j}(z) &= [\alpha_{j-1} -\delta\alpha_{j-1}] z^{j-1}\, 
\ol{\wti E_{j-1} (1/\bar z)} \lb{5.25} 
\end{align} 

By \eqref{2.16} and \eqref{5.7}, 
\begin{equation} \lb{5.26} 
\abs{F_{1,j}(z)} = \wti O(R^{-3j} R^{-j} \max(\abs{z}R^{-1}, 1)^j)  
\end{equation} 
so if $R^3/ (1+\delta)<\abs{z} <R^5$, 
\begin{equation} \lb{5.27} 
\sum_{j=1}^\infty\, \abs{F_{1,j}(z)} \leq \sum_{j=1}^\infty R^{-5j}  
\abs{z}^j <\infty 
\end{equation} 

By \eqref{5.14}, if $\abs{z} >R^3/(1+\delta)$, 
\begin{equation} \lb{5.28} 
z^k \wti E_k (1/\bar z) =\wti O(R^{-n} \max (1,\abs{z} R^{-3})^n)
\end{equation} 
and thus, 
\begin{align} 
\abs{F_{2,j}(z)} &\leq \wti O(R^{-2j} \max (1,\abs{z}R^{-3})^j) \lb{5.29} \\
&\leq \begin{cases} \wti O(\abs{z}^j R^{-5j}) & \text{if } \abs{z} \geq R^3 \\
\wti O(R^{-2j}) & \text{if } \abs{z} \leq R^3 
\end{cases} \lb{5.30} 
\end{align} 
so $\sum_{j=1}^\infty \abs{F_{2,j}(z)}<\infty$ uniformly on compacts of $\{z\mid 
R^3/(1+\delta)\leq \abs{z} <R^5\}$. This implies $r(z)-S(z)-q_3(z)$ is analytic 
there. 
\end{proof} 

\begin{proof}[Proof of Theorem~\ref{T5.1}] As already noted, Theorem~\ref{T2.1} 
proves the results in $\{z\mid 1<\abs{z} <R^3\}$. In $\{z\mid R^3/(1+\delta)  
< \abs{z} <R^5\}$, Theorem~\ref{T5.4} shows $r-S$ is meromorphic with poles 
contained in $\bbG^{(3)} (\{\mu_k\}_{k=1}^K)$. The explicit formula shows that 
there is a pole if there is a single summand contributing to the potential pole. 
\end{proof}

\section{The $R^{2\ell -1}$ Result} \lb{s6} 

In this section, we will prove the following (again, for simplicity of exposition, we 
replace general $P_j(n)$ by constants), which clearly implies Theorems~\ref{T4.6} and 
\ref{T4.8}. 

\begin{theorem}\lb{T6.1} Let $\{\mu_k\}_{k=1}^K$ obey $R\leq \abs{\mu_k} <R^{2\ell-1}$ 
with $\min_k \abs{\mu_k} =R$. Suppose that 
\begin{equation} \lb{6.1} 
\alpha_n = \sum_{k=1}^K c_k \mu_k^{-n-1} + \wti O(R^{-(2\ell-1)})  
\end{equation} 
Then $D(z)^{-1}$ is meromorphic in $\{z\mid \abs{z} < R^{2\ell-1}\}$ with poles 
contained in $\bbG(\{\mu_k\}_{k=1}^K)$. In addition, $S(z) -r(z)$ is meromorphic in 
$\{z\mid 1-\delta_0 < \abs{z} < R^{2\ell+1}\}$ and the only poles in $\{z\mid 
R^{2\ell-1}\leq \abs{z} < R^{2\ell+1}\}$ lie in $\bbG_3 (\{\mu_k\}_{k=1}^K)$. If 
$z_0$ obeys $z_0 =\mu_{i_1}^2 \bar\mu_{i_2}$ with $\abs{z_0} < R^{2\ell+1}$ and 
$z_0$ cannot be written as any other $\bbG(\{\mu_k\}_{k=1}^K)$ product, then $s(z) 
-r(z)$ has a pole at $z_0$. 
\end{theorem} 

The strategy is the same as in the last section. Pick $\delta >0$ so that $\abs{\mu_k} 
< R^{2\ell -1}/(1+\delta)$ for all $k$. We will prove the following estimate on the 
$\Phi$'s and $\Phi^*$'s inductively: 

\begin{theorem} \lb{T6.2} Under the hypothesis of Theorem~\ref{T6.1}, in $\calQ\equiv 
\{z\mid\abs{z} \leq R^{-(2\ell -1)} (1+\delta)\}$, we have 
\begin{equation} \lb{6.2} 
\Phi_n(z) = \sum_p f_{p,\ell}^{(\ell)}(z) \bar w_{p,\ell}^{-n} + 
E_{n,\ell}(z)
\end{equation} 
where the sum is over all points $w$ in 
\[
\biggl[\, \bigcup_{m=1}^{2\ell-3} \bbG^{(2m-1)} (\{\mu_k\}_{k=1}^K)\biggl] 
\bigcap \biggl\{ z\biggm| R\leq \abs{z} < \f{R^{2\ell-1}}{1+\delta}\biggr\} 
\]
each $f_{p,\ell}^{(\ell)}$ is analytic in $\calQ$, and on $\calQ$, 
\begin{equation} \lb{6.3} 
\abs{E_{n,\ell}(z)} = \wti O(\max (R^{-(2\ell-1)},\abs{z})^n) 
\end{equation} 
In addition, in $\calQ$, 
\begin{equation} \lb{6.4} 
\Phi_n^*(z) -d(z)^{-1} = \sum_p g_{p,\ell}(z) y_{p,\ell}^{-n} + 
\wti E_{n,\ell}(z)
\end{equation} 
where the sum is over all products, $p=\mu_{i_1} \dots \mu_{i_m} 
\bar\mu_{i_{m+1}} \dots \bar \mu_{i_{2n}}$ with $i_1, \dots, i_{2m} \in 
\{1, \dots, K\}$ and with the product lying in $\{z\mid 2R \leq \abs{z} 
< R^{2\ell}/(1+\delta)\}$ and on $\calQ$, each $g_{p,\ell}$ is analytic 
and 
\begin{equation} \lb{6.5} 
\abs{\wti E_{n,\ell}(z)} = \wti O(R^{-n} (\max (R^{-(2\ell -1)},z))^n) 
\end{equation} 
\end{theorem} 

\begin{proof} The proof is by induction in $\ell$. Theorems~\ref{T5.2} and \ref{T5.3} 
establish the case $\ell =2$. Suppose that we have the result for $\ell-1$ with 
$\ell \geq 3$ and that \eqref{6.1} holds. Write $\alpha_n$ as in \eqref{5.6} 
where now 
\begin{equation} \lb{6.6} 
(\delta\alpha)_n = \wti O(R^{-(2\ell -1)n}) 
\end{equation} 

In \eqref{5.5}, do the following: 
\begin{align} 
E_{0,\ell}^{(n)} &= z^n \lb{6.7} \\ 
E_{1,\ell}^{(n)} &= -\sum_{j=1}^n \bar\alpha_{n-j} z^{j-1} 
\biggl[\Phi_{n-j}^*- d(z)^{-1} - \sum_p g_{p,\ell-1}(z) y_{p,\ell-1}^{-(n-j)}\biggr] \lb{6.8} \\ 
E_{2,\ell}^{(n)} &= -\sum_{j=1}^n \ol{(\delta\alpha)}_{n-j} z^{j-1} \biggl[ d(z)^{-1} 
+ \sum_p g_{p,\ell-1}(z) y_{p,\ell-1}^{(n-j)}\biggr] \lb{6.9} 
\end{align} 
Then 
\begin{equation} \lb{6.10} 
\Phi_n - \sum_{k=0}^2 E_{k,\ell}^{(n)} = -\sum_{j=1}^n \biggl(\, \sum_{k=1}^K 
\bar c_k \bar\mu_k^{-(n-j)-1}\biggr)\biggl[ d(z)^{-1} + \sum_p g_{p,\ell-1}(z) 
y_{p,\ell-1}^{-(n-j)}\biggr]
\end{equation} 
Define $E_{3,\ell}^{(n)}$ to be the summand on the right side of \eqref{6.10} with 
the sum from $n+1$ to $\infty$. 

The infinite sum yields geometric series which precisely have the form $\sum_p 
g_{p,\ell}(z) y_{p,\ell}^{-n}$, and as in the proof of Theorem~\ref{T5.1}, 
\[
\abs{E_{0,\ell}^{(n)}} + \abs{E_{1,\ell}^{(n)}} + \abs{E_{2,\ell}^{(n)}} + 
\abs{E_{3,\ell}^{(n)}} =\wti O(\max (\abs{z},R^{-(2\ell-1)}))
\]
since $\wti E_{n,\ell-1}$ by induction has $\wti O(\max(R^{-2\ell}, \abs{z})^n)$ 
decay and other terms are bounded by $\wti O(\max (R^{-1}, \abs{z})^n)$. This proves 
\eqref{6.3} for $\ell$. 

To bound $\Phi_n^*(z) -d(z)^{-1}$, we use \eqref{5.15}, replace $\Phi_j(z)$ by 
\eqref{6.2}, $\alpha_n$ by $(\delta\alpha)_n$ plus the asymptotic exponentials, 
and obtain \eqref{6.4} and \eqref{6.5} for $\ell$ by the same estimate as in 
Theorem~\ref{T5.3}. 
\end{proof} 

Basically, using $\Phi_n^*$ to order $R^{-2(\ell-1)n}$ in the expansion of $\Phi_n$ gets 
us $\Phi_n$ to order $R^{-(2\ell-1)n}$, and then plugging that into the expansion of 
$\Phi_n^*$ gets us $\Phi_n^*$ to order $R^{-2\ell n}$. Each full iteration improves 
by $R^{-2n}$. 

\begin{proof}[Proof of Theorem~\ref{T6.1}] By induction, $S-r$ is meromorphic in $\{z 
\mid\abs{z} < R^{2\ell-1}\}$, so knowing $S$ meromorphic implies meromorphicity of $r$  
and so $D^{-1}$ there, and the poles of both $S-r$ and $S$ lie in $\bbG(\{\mu_k\}_{k=1}^K)$. 

Using \eqref{6.4} in $\abs{z} < R^{-(2\ell-1)}/(1+\delta)$ yields an expansion of $\Phi_n(z)$ 
in $\abs{z} > R^{2\ell-1}/(1+\delta)$. Plug this into \eqref{5.21} and use \eqref{5.6}. The 
purely geometric terms sum to poles in $\{z\mid R^{2\ell-1} < \abs{z} < R^{2\ell+1}\}$. The 
bounds on the errors as in the proof of Theorem~\ref{T5.4} converge to an analytic 
function in the annulus. The poles are clearly in $\bbG_3 (\{\mu_k\}_{k=1}^K)$. 

Tracking the contribution of a single $\mu_1^2 \bar\mu_2$ shows that it yields a 
nonvanishing pole which, by the unique product hypothesis, cannot be cancelled. 
\end{proof}

\bigskip

\end{document}